\documentclass[12pt,a4paper]{article}
\usepackage{defns}
\usepackage{graphicx}
\usepackage[colour]{ajr}

\newcommand{\uj}{u_j}
\newcommand{\up}{u_{j+1}}
\newcommand{\um}{u_{j-1}}
\newcommand{\xj}{x_j}
\newcommand{\ibc}{\textsc{ibc}}
\newcommand{\she}{Swift-Ho\-hen\-berg equation}
\newcommand{\gle}{Ginz\-burg-Lan\-dau equation}
\newcommand{\cc}[1]{{\overline #1}}
\newcommand{\cM}{{\cal M}}

\begin{document}

\title{Holistically discretise the \she\ on a scale larger than its 
spatial pattern} 

\author{A. J. Roberts\thanks{Department of Mathematics and Computing, 
University of Southern Queensland, Toowoomba, Queensland 4352, 
\textsc{Australia}.  \protect\url{mailto:aroberts@usq.edu.au}}} 

\maketitle

\begin{abstract}
	I introduce an innovative methodology for deriving numerical 
	models of systems of partial differential equations which exhibit 
	the evolution of spatial patterns.  The new approach directly 
	produces a discretisation for the evolution of the pattern 
	amplitude, has the rigorous support of centre manifold theory at 
	finite grid size~$h$, and naturally incorporates physical 
	boundaries.  The results presented here for the \she\ suggest the 
	approach will form a powerful method in computationally exploring 
	pattern selection in general.  With the aid of computer algebra, the 
	techniques may be applied to a wide variety of equations to derive 
	numerical models that accurately and stably capture the dynamics 
	including the influence of possibly forced boundaries.
\end{abstract}



\tableofcontents

\section{Introduction}

The evolution of spatial patterns is an important area of research.
One example of interest itself and serving as a model for other 
systems is Rayleigh-B\'enard convection 
\cite[e.g.]{Segel69,Greenside84,Newell93} in which a fluid layer is 
heated from below and cooled from above.
The flow flow self organises into an evolving pattern of upwelling and 
downwelling.
One of the most useful models describing such pattern evolution, for 
example see~\cite[(9)]{Cross82}, \cite[(2.11)]{Newell93} 
or~\cite[(13)]{Kramer95}, is the \gle\ which we consider in in one 
spatial dimension:
\begin{equation}
	a_t=ra+ca_{xx}-d|a|^2a
	\label{eq:gle}
\end{equation}
where $a(x,t)$ is the (complex) amplitude of the pattern in any locale, 
subscripts denote partial derivatives, and $r$, $c$ and~$d$ are 
specific constants.
The \gle\ describes the evolution of the complex amplitude of spatially 
periodic structures as they evolve and interact.
The derivation of the \gle\ in any specific pattern problem is 
fundamentally based upon  the underlying 
structure varying slowly in space-time.
However, here in~\S\ref{Scmth} we derive the discrete form
\begin{equation}
	\dot a_j=ra_j +\frac{c}{h^2}(a_{j+1}-2a_j+a_{j-1})-d|a_j|^2a_j
	\label{eq:dgle}
\end{equation}
of the \gle~(\ref{eq:gle}) without ever invoking slow space-time 
variations.
Following earlier research introducing holistic 
discretisation~\cite{Roberts98a}, the derivation is rigorously based 
upon centre manifold theory~\cite[e.g.]{Carr81,Carr83b} and ensures 
the discretisation faithfully models the underlying system.

One amazing consequence of this approach is that we straightforwardly 
incorporate physical boundaries into the discrete model.
Boundaries have significant effects on the pattern evolution 
\cite[e.g.]{Segel69,Cross83}.
Previously there have only been limited, usually just linear 
\cite[p937--8, e.g.]{Cross82}, arguments about how to incorporate such 
boundaries into the analysis to give boundary conditions for the 
\gle~(\ref{eq:gle}).
But here we use nonlinear analyses to derive appropriately modified 
modifications of~(\ref{eq:dgle}) near the boundary to account for its 
effects.
In~\S\ref{Sbc} we discuss time dependent Dirichlet and Neumann 
boundaries as two examples.
Our approach resolves the subgrid fields within each element.
Near a boundary the method resolves the modifications to the subgrid 
fields due to the influence of the boundary and thus naturally creates 
a discretisation appropriate to the specified boundary condition.

As a simple prototype pattern evolution problem, we here consider the 
\she~\cite{Swift77} in one spatial dimension:
\begin{equation}
	u_t=ru-(1+\partial_x^2)^2u-u^3\,.
	\label{eq:she}
\end{equation}
This equation has often been used, as we do here, to investigate 
issues in pattern selection \cite[e.g.]{Cross84, Cross86}, especially 
the influence of physical boundaries in 2D \cite[e.g.]{Greenside84, 
Bestehorn89}.
Linearly, the \she~(\ref{eq:she}) has spatially periodic solutions 
$u\propto \exp( \lambda t+ikx)$ for wavenumber~$k$ where the 
growth-rate
\begin{equation}
	\lambda=r-(1-k^2)^2\,.
	\label{eq:shlam}
\end{equation}
Thus for parameter~$r<0$ all spatial modes decay, whereas for $r>0$ a 
band of modes near $|k|=1$ may grow.
Exactly at the critical parameter~$r=0$, all modes decay except for 
the two neutral modes $\exp(\pm ix)$.
We construct a holistic numerical discretisation based upon these 
neutral modes.

\section{Develop a centre manifold discretisation}
\label{Scmth}

Here we develop a model of the system using invariant manifold theory.
The analysis is essentially local in space, it resolves structures 
over a number of ``convective rolls''.
Previous use of invariant manifold theory~\cite{Graham93} has 
generally sought global models evolving on the so-called inertial 
manifold~\cite[e.g.]{Temam90}.

\begin{figure}[tbp]
	\centering
	\includegraphics{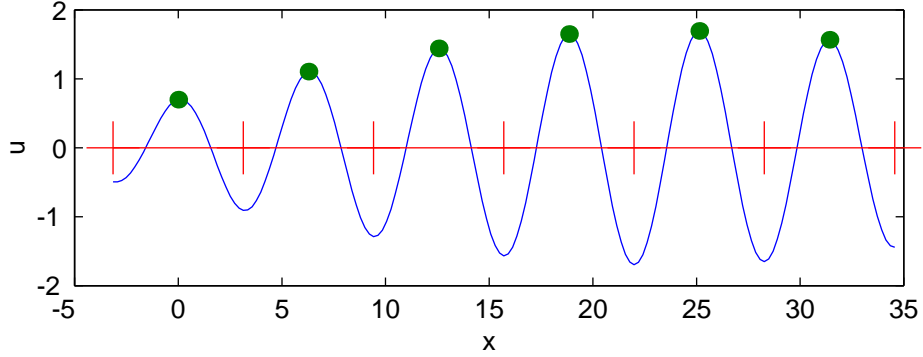}
	\caption{schematics diagram showing a varying ``roll'' structure 
	(solid curve) discretised by introducing artificial internal 
	boundaries every period ($p=1$) at odd multiples of~$\pi$ 
	(vertical bars).  The ``roll'' field in the $j$th~element is 
	parametrised by the amplitude~$a_j$ (discs).}
	\label{fig:rolls}
\end{figure}

Exactly at criticality, long lasting solutions are $2\pi$-periodic in~$x$.
As shown schematically in Figure~\ref{fig:rolls}, introduce artificial 
internal boundaries every $p$~periods, spacing $h=p2\pi$.
These divide the domain into elements, each centred on a grid 
point~$\xj$.
Denote the field in the $j$th~element by $\uj(x,t)$.
The internal boundary conditions (\ibc{}s) to hold at each end of an 
element are chosen to be the nonlocal conditions
\begin{eqnarray}
	\uj+\D x\uj & = & (1-\gamma)\left[\uj+\D x\uj\right]_{x=\xj-h/2} 
	\nonumber\\&&{}
	+\gamma\left(\up+\D x\up\right)\quad \mbox{at $x=\xj+h/2$\,,}
	\label{eq:rbc}  \\
	\uj-\D x\uj & = & (1-\gamma)\left[\uj-\D x\uj\right]_{x=\xj+h/2} 
	\nonumber\\&&{}
	+\gamma\left(\um-\D x\um\right)\quad \mbox{at $x=\xj-h/2$\,,}
	\label{eq:lbc}
\end{eqnarray}
and the same for the second derivative $v=u_{xx}$\,.
These \ibc{}s are parametrised by~$\gamma$ which controls the 
interaction and information flow between adjacent elements.
When this coupling parameter~$\gamma=0$ these \ibc{}s reduce to 
conditions,
\begin{equation}
	\left[\uj\pm\D x\uj\right]_{\xj-h/2}
	=\left[\uj\pm\D x\uj\right]_{\xj+h/2}\,,
	\label{eq:pbc}
\end{equation}
requiring the solution in the $j$th~element is $h$-periodic---the 
field and its derivatives at the right end of the element must match 
smoothly the field at the left end.
Whereas when $\gamma=1$ these \ibc{}s require continuity of the field 
and its derivatives between abutting ends of adjacent elements,
\begin{displaymath}
	\uj\pm\D x\uj=\up\pm\D x\up\quad\mbox{at $x=\xj+h/2$\,.}
\end{displaymath}
Thus when the coupling parameter~$\gamma=0$ each element is isolated 
from all others, whereas when $\gamma=1$ we ensure enough continuity 
between elements to recover the \she~(\ref{eq:she}) throughout the 
domain.
We use these \ibc{}s to model the \she, when~$\gamma=1$, by basing 
analysis on the discrete elements that are isolated at~$\gamma=0$.

Apply centre manifold theory based upon the linear dynamics in 
uncoupled elements.
By notionally adjoining the dynamically trivial equations 
$\dot\gamma=\dot r=0$ we treat all terms with a factor of the coupling 
parameter~$\gamma$ or the forcing parameter~$r$ as ``nonlinear'' 
perturbations.
$h$-periodic solutions of the linearised \she, 
$u_t=-(1+\partial_x^2)^2u$, then are $\exp(\lambda_nt+ik_nx)$ for any 
integer~$n$ with wavenumber~$k_n=2\pi n/h=n/p$ and 
growth-rate~$\lambda_n=-(1-k_n^2)^2$\,.
Thus within each element all modes decay exponentially except the two 
modes with wavenumber~$k_{\pm p}=\pm 1$ which are neutral.
Hence ``linearly,'' the solution in each element evolves exponentially 
quickly in time to the roll solution
\begin{equation}
	u_j=a_je^{ix}+b_je^{-ix}\,,
	\label{eq:lins}
\end{equation}
where $a_j$~and~$b_j$ are the complex amplitudes of the rolls---the 
solution is real if and only if all~$b_j$ are the complex conjugate 
of~$a_j$.
Theory \cite[e.g.]{Carr81,Carr83b} then assures us that for the 
original nonlinear \she, coupled between elements by the 
\ibc{}s~(\ref{eq:rbc}--\ref{eq:lbc}), there exists a centre 
manifold~$\cM$ on which solutions are parametrised by the collection 
of evolving amplitudes~$\vec a$ and~$\vec b$:
\begin{eqnarray}
	&& 
	u=u_j(\vec a,\vec b,\gamma,r,x) 
    \label{eq:field}\\
	\mbox{such that}&&
	\dot a_j=g_j(\vec a,\vec b,\gamma,r)
	\quad\mbox{and}\quad
	\dot b_j=\bar g_j(\vec a,\vec b,\gamma,r)\,.
	\label{eq:ansatz}
\end{eqnarray}
Moreover, theory assures us quite generally that solutions of the 
\she\ exponentially quickly approach solutions of~(\ref{eq:ansatz}).
Lastly, theory asserts that the functions 
in~(\ref{eq:field}--\ref{eq:ansatz}) are determined by substitution 
into the \she\ and the internal boundary conditions and then solving 
to some order in the parameters.

Computer algebra available from the author performs the tedious 
calculations using an iterative algorithm~\cite{Roberts96a}.
However, before undertaking the modelling we define the amplitudes 
precisely as the element averages
\begin{equation}
	a_j(t)=\frac{1}{h}\int_{x_j-h/2}^{x_j+h/2} u(x,t)e^{-ix}dx
	\quad\mbox{and}\quad
	b_j(t)=\frac{1}{h}\int_{x_j-h/2}^{x_j+h/2} u(x,t)e^{+ix}dx\,.
	\label{eq:ampl}
\end{equation}
Then the solution field is found to be
\begin{eqnarray}
	u_j&=&\phantom{{}+}a_je^{+ix}
	+\frac{\gamma}{4h}e^{+ix}\left[ \left( \delta^2a_j-2i\mu\delta b_j 
	\right) +\left( 4\mu\delta a_j -2i\delta^2b_j \right)x\right]
	\nonumber\\&&{}
	+b_je^{-ix}
	+\frac{\gamma}{4h}e^{-ix}\left[ \left( \delta^2b_j+2i\mu\delta a_j 
	\right) +\left( 4\mu\delta b_j +2i\delta^2a_j \right)x\right]
	\nonumber\\&&{}
	+\Ord{\gamma^3+A^3+r^{3/2}}\,,
	\label{eq:ufld}
\end{eqnarray}
where $A=\|\vec a\|+\|\vec b\|$ measures the overall amplitude of the 
solution field\footnote{A multinomial term $\gamma^lA^mr^n$ is 
$\Ord{\gamma^q+A^q+r^{q/2}}$ if and only if $l+m+2n\geq q$\,.
} and where $\mu$~and~$\delta$ are central mean and difference 
operators, $\mu a_j =(a_{j+1/2} +a_{j-1/2})/2$ and $\delta a_j 
=a_{j+1/2} -a_{j-1/2}$ respectively.
These $\Ord{\gamma}$ terms show the leading order effect of 
neighbouring elements upon the subgrid field in each element; however, 
there is as yet no effect upon the evolution which is still neutral 
$g_j=\bar g_j=\Ord{\gamma^2}$.
Using computer algebra to compute the next order of interactions and I 
find the evolution
\begin{eqnarray}
\dot a_j & = & ra_j +\frac{4\gamma^2}{h^2}\delta^2a_j
-3\gamma^2a_j^2b_j +\Ord{\gamma^4+A^4+r^2}\,, \label{eq:aev} \\
	\dot b_j & = &  rb_j +\frac{4\gamma^2}{h^2}\delta^2b_j
	-3\gamma^2a_jb_j^2 +\Ord{\gamma^4+A^4+r^2}\,.
	\label{eq:bev}
\end{eqnarray}
Evaluating these at the coupling parameter~$\gamma=1$ we recover a 
numerical model for the \she\ which is just the discrete 
version~(\ref{eq:dgle}) of the appropriate \gle~(\ref{eq:gle}).

In this approach there is little merit in discussing equivalent 
\pde{}s to the derived discrete models~(\ref{eq:aev}--\ref{eq:bev}).
The reason is that here the element size~$h$ must contain an integral 
number of rolls and so the limit~$h\to0$ is not physically valid.
For the same reason, there appears to be no merit in seeking 
consistency with a \pde\ as done for other holistic 
discretisations~\cite{Roberts00a}.
In this approach we should discuss the numerical model as it stands.
But it is relevant to observe that the long-wave limit, slow 
variations in~$j$, should and does reduce to the relevant 
\gle~(\ref{eq:gle}).

\section{Boundary conditions are straightforwardly determined}
\label{Sbc}

Now consider the boundaries to the physical domain.
For simplicity suppose that the conditions at the boundary, say the 
left boundary, are either one of the two cases:
\begin{eqnarray}
	u=(-1)^p\alpha(t) &\mbox{and}& u_{xx}=(-1)^p\beta(t)\,;
	\label{eq:evenbc}\\
	\mbox{or}\quad
	u_x=(-1)^p\alpha(t) &\mbox{and}& u_{xxx}=(-1)^p\beta(t)\,.
	\label{eq:oddbc}	
\end{eqnarray}
These physical boundary conditions are incorporated into the analysis 
by replacing the left-hand \ibc~(\ref{eq:lbc}) of the left-most 
element, say element~$j=1$, by a boundary condition corresponding to 
one of~(\ref{eq:evenbc}) or~(\ref{eq:oddbc}).
Consequently this left physical boundary is, without loss of 
generality, located at $x=x_1-h/2$\,.
However, we want this left-most element, $j=1$, to still have the 
neutral periodic solution~(\ref{eq:lins}) when the inter-element 
coupling parameter~$\gamma=0$.
Thus we actually replace the left-hand \ibc~(\ref{eq:lbc}) of the 
left-most element by
\begin{eqnarray}
	u_1-\D x{u_1} &=& (1-\gamma)\left[{u_1}-\D x{u_1}\right]_{x=x_1+h/2} 
	\nonumber\\&&{}
	\mp\gamma\left({u_1}+\D x{u_1}\right) \pm2(-1)^p\gamma\alpha(t)
	\quad\mbox{at $x=x_1-h/2$,}
	\label{eq:lmbc}
\end{eqnarray}
and similarly for~$v=u_{xx}$.
When the coupling parameter~$\gamma=0$ (\ref{eq:lmbc}) reduces to 
requiring $h$-periodic solutions.
Whereas when $\gamma=1$ (\ref{eq:lmbc}) reduces to
\begin{displaymath}
	(1\pm1)u_1-(1\mp1)\D x{u_1}=\pm2(-1)^p\alpha
	\quad\mbox{at $x=x_1-h/2$}\,,
\end{displaymath}
and similarly for~$v=u_{xx}$, which for the upper choice of signs 
specifies the even derivative \textsc{bc}~(\ref{eq:evenbc}) and for 
the lower choice of signs specifies the odd derivative 
\textsc{bc}~(\ref{eq:oddbc}).
Other specific boundary conditions may be treated similarly, but here 
I just restrict attention to these two cases.

Computer algebra constructs the centre manifold model for both 
alternative boundary conditions: (\ref{eq:evenbc}) giving the upper 
alternative of plus\slash minus signs seen below; and~(\ref{eq:oddbc}) 
giving the lower alternative.
The leading order influence of the boundary forcing upon the field in 
the leftmost element is
\begin{eqnarray}
	u_1 & = & e^{ix}a_1
	+\frac{\gamma}{4h}e^{ix}\left[ (-(2\pm i)a_1+a_2\mp b_1-ib_2) 
	\right.\nonumber\\&&\qquad\left.{}
	+2(\pm ia_1+a_2\pm(1\pm 2i)b_1-ib_2)x \right]
	\nonumber\\&&{}
	\pm\frac{\gamma^2\alpha}{h}e^{ix}\left[ \frac{7+5i}{16} 
	-\frac{2+3i}{4}x +\frac{1-i}{96}(h^2-12x^2) \right]
	\nonumber\\&&{}
	\pm\frac{\gamma^2\beta}{h}e^{ix}\left[ \frac{3+i}{16} 
	-\frac{i}{4}x +\frac{1-i}{96}(h^2-12x^2) \right]
	\nonumber\\&&{}
	+e^{-ix}b_1
	+\frac{\gamma}{4h}e^{-ix}\left[ (\mp a_1+ia_2-(2\mp i)b_1+b_2) 
	\right.\nonumber\\&&\qquad\left.{}
	+2(\pm(1\mp 2i)a_1+ia_2\mp ib_1+b_2)x \right]
	\nonumber\\&&{}
	\pm\frac{\gamma^2\alpha}{h}e^{-ix}\left[ \frac{7-5i}{16} 
	-\frac{2-3i}{4}x +\frac{1+i}{96}(h^2-12x^2) \right]
	\nonumber\\&&{}
	\pm\frac{\gamma^2\beta}{h}e^{-ix}\left[ \frac{3-i}{16} 
	+\frac{i}{4}x +\frac{1+i}{96}(h^2-12x^2) \right]
	\nonumber\\&&{}
	+\Ord{\gamma^3+A^3+r^{3/2},|\dot\alpha|+|\dot\beta|}
	\label{eq:fld1}
\end{eqnarray}
\begin{figure}[tbp]
	\centering
	\includegraphics{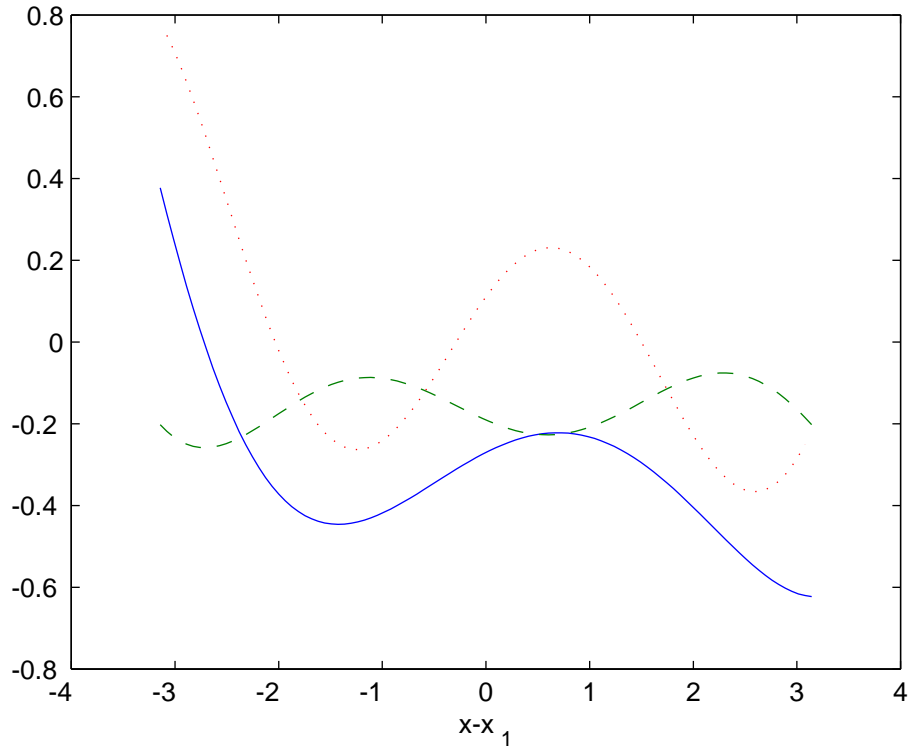}
	\caption{spatial subgrid structure~(\ref{eq:fld1}) of the 
	coefficients of~$\alpha$ (solid), $\beta$ (dashed) and its second 
	derivative (dotted) when the boundary condition at~$x=x_1-\pi$ 
	is~(\ref{eq:evenbc}) for an element size~$h=2\pi$\,.  These give 
	the field if all amplitudes~$a_j=b_j=0$\,.  See the $\alpha$ and 
	the second derivative curves are at the boundary~$x=0$ about a 
	value of one higher than the field in the bulk of the 
	element---this is appropriate as~$\alpha$ specifies the boundary 
	value and $\beta$ the second derivative.  }
	\label{fig:bcelem}
\end{figure}%
See there are two sorts of effects of the boundary: first, the subgrid 
field generated by the grid values of the leftmost two elements 
differs from that in an interior element~(\ref{eq:ufld}); secondly, 
the forcing parametrised by $\alpha$~and~$\beta$ also adjusts the 
subgrid field.
For example, in Figure~\ref{fig:bcelem} is plotted the basis of the 
field in the leftmost element when all grid values~$a_j=b_j=0$\,.
See that they resolve a boundary layer structure on a scale of $\Delta 
x\approx 1$---the resolution is crude because this is only the first 
approximation to the adjustment necessary to account for the boundary.
By resolving the subgrid spatial structure near a boundary we will 
generate appropriate discretisations for the given boundary condition.
The resolution of structure near the boundary is important for 
deriving correct boundary conditions for models~\cite[p937, 
e.g.]{Cross82}.
As demonstrated in Figure~\ref{fig:bcelem}, our approach does some of 
the near boundary analysis identified in~\cite{Roberts92c} as 
necessary for deriving boundary conditions for mathematical models 
such as the differential \gle.

The relevant evolution is obtained by analysing to quadratic terms in 
the coupling parameter~$\gamma$.
We find that $\dot a_2$ and~$\dot b_2$ are identical to that obtained 
for the interior, (\ref{eq:aev}--\ref{eq:bev}) with $j=2$.
However, the evolution for the amplitude of the rolls in the leftmost 
element adjacent to the boundary is
\begin{eqnarray}
\dot a_1 & = & ra_1 +\frac{4\gamma^2}{h^2}(a_2-2a_1\mp b_1) -3a_1^2b_1 
      \mp\frac{\gamma^2}{h}(1-i)(\alpha+\beta)
\nonumber\\&&{}
+\Ord{\gamma^4+A^4+r^2,|\ddot\alpha|+|\ddot\beta|}\,, \label{eq:a1ev} \\
	\dot b_1 & = &  rb_1 +\frac{4\gamma^2}{h^2}(b_2-2b_1\mp a_1) -3a_1b_1^2 
    \mp\frac{\gamma^2}{h}(1+i)(\alpha+\beta)
\nonumber\\&&{}
	+\Ord{\gamma^4+A^4+r^2,|\ddot\alpha|+|\ddot\beta|}\,.
	\label{eq:b1ev}
\end{eqnarray}
Set the inter-element coupling parameter~$\gamma=1$ to recover a 
discretisation.
In the absence of boundary forcing, $\alpha=\beta=0$, we analyse a 
little of the dynamics near the boundary by supposing for illustrative 
purposes that $a_1=a_2=a$ and $b_1=b_2=\cc a$ (for real solutions).
The evolution equations~(\ref{eq:a1ev}--\ref{eq:b1ev}) then reduce to
\begin{equation}
	\dot a\approx ra-\frac{4}{h^2}(a\pm \cc a)-3|a|^2a\,.
	\label{eq:abdry}
\end{equation}
First, consider the upper alternative when $u$~and~$u_{xx}$ are 
specified at the boundary, (\ref{eq:evenbc}).
Linearly, the real part of~$a$ will decay at a rapid rate, a negative 
growth-rate of~$r-8/h^2$; whereas the imaginary part has a 
growth-rate~$r$.
Thus the amplitudes near the boundary will rapidly evolve to be pure 
imaginary.
This corresponds to solutions~$u\propto\sin(x)$ as expected for the 
given boundary conditions.
Secondly, consider the converse lower alternative when 
$u_x$~and~$u_{xxx}$ are specified at the boundary, (\ref{eq:oddbc}).
Linearly, the imaginary part of~$a$ will decay at a rapid rate, a 
negative growth-rate of~$r-8/h^2$; whereas the real part has a 
growth-rate~$r$.
Thus the amplitudes near the boundary will rapidly evolve to be purely 
real.
This corresponds to solutions~$u\propto\cos(x)$ as expected for the 
given boundary conditions.
Thus the near boundary discretisation~(\ref{eq:a1ev}--\ref{eq:b1ev}) 
selects for fields~$u$ with evenly spaced rolls located so that~$u$ is 
zero at the boundary in the case of specified even derivatives, and so 
that the derivatives of~$u$ are zero at the boundary in the case of 
specified odd derivatives.
These properties agree with earlier more specific work 
\cite[e.g.]{Segel69}

The presence of boundary forcing will push the system away from its 
otherwise natural equilibria.
For example, for the upper alternative, (\ref{eq:abdry}) breaks the 
symmetry in $\Im(a)$ and predicts an equilibrium $\Re(a)\approx 
-h(\alpha+\beta)/8$ instead of zero and so there is a change in 
amplitude and phase of the rolls near the boundary.
Observe that time variations in the forcing, $\dot\alpha$ 
and~$\dot\beta$, are not significant at this level of approximation: 
this absence is surprising given the appearance of time derivatives of 
forcing in other boundary discretisations at finite grid 
size~\cite{Roberts01b}; the absence is linked to the 
definition~(\ref{eq:ampl}) of the amplitudes because other amplitude 
definitions generate $\dot\alpha$~and~$\dot\beta$ terms in the model.
Cox and Roberts~\cite{Cox91} showed that special parametrisations 
could eliminate time dependence upon forcing in a centre 
manifold---evidently the amplitude definition~(\ref{eq:ampl}) at least 
approximates such a parametrisation.
The important point is that our treatment of near boundary elements 
generates discretisations which naturally incorporate the forced 
boundary conditions.

\section{Conclusion}

We have introduced a rigorously based method for deriving discrete 
amplitude equations for the evolution of spatial patterns directly 
from the original \pde{}s.
The method has been applied the the \she\ but the principles apply to 
a wide variety of \pde{}s governing spatio-temporal evolution.
One remarkable advantage of this direct derivation is that we may also 
derive appropriate discretisations near any forced boundary.

The generalisation to patterns in 2D space appears a straightforward 
generalisation of that employed for reaction-diffusion 
equations~\cite{MacKenzie00b}.
Space would be tessellated in some regular manner by introducing 
appropriately periodic \ibc{}s.
Then a centre manifold model would be constructed with a finite number 
of basis wavevectors in each element analogous to the regular analysis 
of Skeldon \etal~\cite{Skeldon98}.

Computer algebra readily determines models of higher-order in the 
asymptotic expansions, higher-order in both or either of nonlinearity 
or of stencil width.
However, the \ibc{}s~(\ref{eq:rbc}--\ref{eq:lbc}) used then here 
appear to break a symmetry of the \she.
More research is needed into a better form of the \ibc{}s as well as 
further applications of the methodology.

\paragraph{Acknowledgement:} this research is partially supported by a 
grant from the Australian Research Council.

\bibliographystyle{plain}
\bibliography{ajr,bib,new}

%
%

\end{document}